\documentclass{article}

  \usepackage{amsmath}
  \usepackage{amsfonts}
  \usepackage{paralist}
  \usepackage{graphics} 
  \usepackage{epsfig} 
  \usepackage{overpic}

 \usepackage[colorlinks=true]{hyperref}

\newtheorem{theorem}{Theorem}[section]
\newtheorem{corollary}{Corollary}
\newtheorem{lemma}[theorem]{Lemma}
\newtheorem{proposition}{Proposition}

\newtheorem{definition}[theorem]{Definition}
\newtheorem{remark}{Remark}

\title{Dynamics of Functions with an Eventual Negative Schwarzian Derivative}

\author{Benjamin Webb}

\title{text}

\begin{document}

\centerline{\large{\textbf{Dynamics of Functions with an Eventual Negative }}}
\centerline{\large{\textbf{Schwarzian Derivative}}}
\bigskip
\bigskip

\centerline{\scshape Benjamin Webb }
\medskip
{\footnotesize
 \centerline{School of Mathematics, Georgia Institute of Technology}
   \centerline{Atlanta, GA 30322, USA}
}

\bigskip


\begin{abstract}
In the study of one dimensional dynamical systems one often assumes that the functions involved have a negative Schwarzian derivative. In this paper we consider a generalization of this condition. Specifically, we consider the interval functions of a real variable having some iterate with a negative Schwarzian derivative and show that many known results generalize to this larger class of functions. The introduction of this class was motivated by some maps arising in neuroscience.
\end{abstract}

\section{Introduction}

Singer was the first to observe that if a function has a negative Schwarzian derivative then this property is preserved under iteration and moreover that this property puts restrictions on the type and number of periodic orbits the function can have \cite{1}. These properties, as well as those results derived from them, essentially rely on the global structure functions with a negative Schwarzian derivative have. For a full list of such properties see \cite{14}.   

Later it was found that such functions also possess local properties useful in establishing certain distortion bounds. For the most part these properties are concerned with the way in which functions with a negative Schwarzian derivative increase cross ratios. Because of these special properties some results are known only in the case where a function has a negative Schwarzian derivative. An exception, however, to this is the result by Kozlovski \cite{2} where it was shown that the assumption of a negative Schwarzian derivative is superfluous in the case of any $C^3$ unimodal map with nonflat critical point. For such functions there is always some interval around their critical value on which the first return map has a negative Schwarzian derivative. That is locally all such $C^3$ maps behave as maps with a negative Schwarzian derivative. More recently using the same technique van Strien and Vargas have generalized this result to multimodal functions \cite{7}. Also Graczyk, Sands and Swiatek have shown that any $C^3$ unimodal map with only repelling periodic points is analytically conjugate to a map with a negative Schwarzian derivative \cite{5}. The main purpose of these results is to relate functions without a negative Schwarzian derivative to functions with this property.

In this paper we consider those $C^2$ functions on a finite interval of the real line having some iterate with a negative Schwarzian derivative. This class, which we call functions with an \textit{eventual negative Schwarzian derivative}, was originally introduced by L. Bunimovich (2007, personal communication) in an attempt to describe some one-dimensional maps which appear in neuroscience \cite{10,11}. It is noteworthy that this class of functions is broader than those previously considered in the study of unimodal and multimodal maps related to functions with a negative Schwarzian derivative (see for example \cite{2,7,5}). To demonstrate this we will present examples of such functions and further examples can also be found in \cite{15}. Moreover, as having an eventual negative Schwarzian derivative is not an asymptotic condition, verifying whether a function has this property can often be done by direct computation. Hence, the concept of an eventual negative Schwarzian derivative has potential to be very useful in applications. 

Since this class contains functions with a negative Schwarzian derivative as a subset we do not attempt to prove stronger results then have already been proved for this smaller class of functions. With this in mind the large majority of the results we will present will simply be restatements of those already known with the modification that only some iterate of our function need have a negative Schwarzian derivative. We note that this by no means is meant to be an exhaustive list of such results as the purpose of this paper is to give further evidence that the useful properties possessed by functions with negative Schwarzian derivative are not limited to this family of functions.

This paper is organized as follows: In the next section we formally introduce the class of functions we are considering and present our main results. Section 3 presents the proof of those results that are topological in nature. Specifically, we prove an analogue of Singer's theorem in \cite{1} and mention some important concepts and corollaries that will be useful in what follows. Section 4 is comprised of those proofs which are more measure theoretic. Specifically, we generalize the main results in \cite{8,4,9} to this larger class of functions. In section 5 we give a partial characterization of functions that have an eventual negative Schwarzian derivative as well as some examples. The next section is devoted to an application of these results to a one-parameter family of maps that model the electrical activity in a neuronal cell near the transition to bursting \cite{10,11}. Section 7 contains some concluding remarks.

\section{Iterates and the Schwarzian Derivative}
  
\noindent In this paper we consider only $C^2$ functions $f:I\rightarrow I$ on some nontrivial compact interval $I$ of real numbers having a nonempty but finite set of critical points $\mathcal{C}(f)$ and use the notation $f^{\prime}$ to denote the derivative with respect to the spacial variable.

\begin{definition} A $C^2$ function $f:I\rightarrow I$ is said to have a \textit{negative Schwarzian derivative} if on any open interval $J\subset I$, not containing critical points of $f$, $|f^{\prime}|^{-1/2}$ is strictly convex on $J$.
\end{definition}

If $f:I\rightarrow I$ is $C^3$ then its \textit{Schwarzian derivative} is defined off of $\mathcal{C}(f)$ by
\begin{equation} 
S(f)(x)=\frac{f^{\prime\prime\prime}(x)}{f^{\prime}(x)}-
\frac{3}{2}\Big(\frac{f^{\prime\prime}(x)}{f^{\prime}(x)}\Big)^2
\label{eq:no8}
\end{equation}
If $S(f)(x)$ is strictly negative on open intervals not containing critical points of $f$ then $|f^\prime(x)|^{-1/2}$ is strictly convex on these sets or $f$ has a negative Schwarzian derivative. However, the converse does not always hold. 

\begin{definition} We say a $C^2$ function $f:I\rightarrow I$ has an \textit{eventual negative Schwarzian derivative} if there exists $k\in\mathbb{N}$ such that $f^k$ has a negative Schwarzian derivative. The smallest such number $k$ is said to be the \textit{order} of the derivative.
\end{definition}

\begin{definition} A $C^2$ map $f:I\rightarrow I$ is called \textit{S-multimodal} if:\\
(i) $\mathcal{C}(f)$ is nonempty and finite.\\
(ii) For every $c\in \mathcal{C}(f)$, $c$ is \textit{nonflat}. That is, there exists some $\ell\in(1,\infty)$ and $L\in(0,\infty)$ such that
\begin{equation}
\lim_{x\rightarrow c} \frac{|f^{\prime}(x)|}{|x-c|^{\ell-1}}=L,
\label{eq:no1}
\end{equation} where $\ell$ is the \textit{order} of the critical point.\\
\noindent(iii) $f$ has a negative Schwarzian derivative on $I$.\\ 
If $f$ has a single critical point $c$ and $f^{\prime\prime}(c)\neq 0$ we call the map \textit{S-unimodal}.
\end{definition}

\begin{definition} A function $f:I\rightarrow I$ is $S^k$\textit{-multimodal} if there exists a smallest $k\in\mathbb{N}$ such that $f^k(x)$ is $S$-multimodal. If $f$ is unimodal and has this property we say it is \textit{$S^k$-unimodal}
\end{definition}

Note that by the terminology above $S$-unimodal functions are special $S$-multimodal functions and $S^k$-unimodal functions are special $S^k$-multimodal functions. 

In what follows we will make use of the following standard terminology. The \textit{basin} of a periodic point $x$ is the set of points whose forward orbits accumulate on $x$, and $x$ is said to be \textit{attracting} if its basin contains an open set. The \textit{immediate basin} of $x$ is the union of connected components of its basin that contain a point of the orbit of $x$. Furthermore, we say that the periodic point $x$ of order $p$ is a \textit{hyperbolic attractor} if $|(f^p)^{\prime}(x)|<1$, a \textit{hyperbolic repeller} if $|(f^p)^{\prime}(x)|>1$, and \textit{neutral} if $|(f^p)^{\prime}(x)|=1$. Note that it is possible for neutral periodic points to be attracting from one or both sides. We now state the main results of this paper.

\begin{theorem}\label{theorem 1}
If $f:I\rightarrow I$ is $S^k$-multimodal, then the immediate basin of attraction of any attracting periodic point contains either a critical point of $f$ or a boundary point of $I$. Furthermore, any neutral periodic point of $f$, except possibly on the boundary of $I$, is attracting and there exist no interval of periodic points.
\end{theorem}

This type of theorem, often called Singer's theorem as it resembles the result in \cite{1}, is proved in \cite{14} under the assumption that $f$ is $C^3$ with $S(f)(x)<0$. In \cite{12} theorem \ref{theorem 1} is proved in the case that $f$ is $S^1$-multimodal.
  
We note here that theorem \ref{theorem 1} indicates a few of the properties that sets this collection of functions apart from those considered elsewhere. In fact the property of having an eventual negative Schwarzian derivative cannot be generalized by looking at first return maps in the sense that it is global in nature and as first return maps generically introduce discontinuities, i.e. more boundary points, this global structure is not preserved. 

The next result is not a generalization of a known result rather it is a corollary to the main result in \cite{8} using the restriction on the periodic orbits obtained in the previous theorem to simplify the hypothesis. 

Denote $D_n(c)=|(f^n)^{\prime}(f(c))|$. Also, for some measure $\mu$ let $\varphi$ and $\psi$ be bounded H\"{o}lder continuous functions on $I$ and denote the \textit{$n^{th}$ correlation function} by
$$C_n=C_n(\varphi,\psi)=\left|\int(\varphi\circ f^n)\psi d\mu-\int \varphi d\mu \int \psi d\mu \right|.$$
 
\begin{theorem}\label{theorem 2} 
Let $f:I\rightarrow I$ be $C^3$ with an eventual negative Schwarzian derivative, and a finite nonempty critical set $\mathcal{C}(f)$. If every point of $\mathcal{C}(f)$ has order $\ell\in (1,\infty)$ and $f$ satisfies 
\begin{equation}
\sum_n D_n^{-1/(2\ell-1)}(c)< \infty \ \text{for all} \ c\in \mathcal{C}(f)
\label{eq:no2} 
\end{equation} then there exists an $f$-invariant probability measure $\mu$  absolutely continuous with respect to Lebesgue measure (an ``\textit{acip}''). Furthermore, some iterate of $f$ is mixing and $C_n$ decays at the following rates:
\begin{center}
\textit{Polynomial}: If there is $C>0$, $\alpha>2\ell-1$ such that $D_n(c)\geq Cn^\alpha$\\ for all $c\in\mathcal{C}(f)$ and $n\geq 1$, then for any $\tilde{\alpha}<\frac{\alpha-1}{\ell-1}-1$ we have $C_n=\mathcal{O}(n^{-\alpha})$.
\end{center}
\begin{center}
\textit{Exponential}: If there is $C>0$, $\beta>0$ such that $D_n(c)\geq Ce^{\beta n}$\\ for all $c\in\mathcal{C}(f)$ and $n\geq 1$, then there is a $\tilde{\beta}>0$ such that $C_n=\mathcal{O}(e^{-\tilde{\beta}})$.
\end{center}
\end{theorem}

In the original result given in \cite{8} the function $f$ was assumed to have no attracting or neutral periodic points. In effect, the previous theorem says that instead of requiring the function to have no such periodic orbits we may assume that it has an eventual negative Schwarzian derivative. The advantage here is that having an eventual negative Schwarzian derivative is a nonasymptotic condition so it is potentially easier to verify a function has this property than to show it has no attracting or neutral periodic points by some other means.

However, if all that is needed is the existence of an acip, we have different orders of critical points, or the function is $C^2$ but not $C^3$ then we may use the following result which generalizes the main result in \cite{4} to functions with an eventual negative Schwarzian derivative.

\begin{theorem}\label{theorem 3} If  $f$ is an $S^k$-multimodal function and satisfies the condition 

\begin{equation}
\sum_n D_n^{-1/(\ell_{max})}(c)< \infty \ \text{for all}  \ c\in \mathcal{C}(f)
\label{eq:no3} 
\end{equation} where $\ell_{max}$ is the largest order of the critical points in $\mathcal{C}(f)$ then $f$ admits an absolutely continuous invariant probability measure.
\end{theorem}

The next theorem deals with one-parameter families of maps. In the theory presented in \cite{9} there is a special class of functions denoted by $\mathcal{M}$ having strongly expansive properties. This set of functions is helpful in proving under what conditions one-parameter families of functions have absolutely continuous invariant measures for positive Lebesgue measure sets of parameters. Specifically, a technical but generically satisfied transversality condition in the parameter which will be denoted by (\textit{PT}) is required for this to be the case. We refer the reader to the article for details.

\begin{theorem}\label{theorem 4}
Let $f_a:I\rightarrow I$ be a one-parameter family of $C^3$ functions where $a$ belongs to some interval $A$ of the real line. If, for some parameter value $\beta\in A$, $f_{\beta}$ has a finite critical set $\mathcal{C}(f_{\beta})$ and\\
(i) $f_{\beta}$ has an eventual negative Schwarzian derivative of order $k$\\
(ii) $f_{\beta}^{\prime\prime}(c)\neq 0$, $c\in \mathcal{C}(f_{\beta})$\\
(iii) if $f_{\beta}^m(x)=x$, then $|(f_{\beta}^m)^{\prime}(x)|>1$\\
(iv) inf$_{i>0}d(f_{\beta}^i(c),\mathcal{C}(f_{\beta}))>0$, $c\in \mathcal{C}(f_{\beta})$\\
then $f_\beta$ has an absolutely continuous invariant measure. In particular $f^k_{\beta}\in\mathcal{M}$ and if $f^k_{\beta}$ satisfies the condition (\textit{PT}) then on a positive Lebesgue measure set of parameters the family of functions $f_{a}$ has an absolutely continuous invariant probability measure.
\end{theorem}

Theorem \ref{theorem 4} is an extension of the results in \cite{9} with the modification that we require only an eventual negative Schwarzian derivative of order $k\geq 1$ instead of a negative Schwarzian derivative.

\section{Topological Properties} In this section we prove theorem \ref{theorem 1} along with some corollaries that will be needed in the following sections but are of interest in their own right. We now give a proof of theorem $\ref{theorem 1}$.\\

\noindent\textit{Proof.}
Let $f:I \rightarrow I$ be an $S^k$-multimodal function of order $k\geq 1$. From the proof of Singer's theorem in \cite{12} the results of theorem \ref{theorem 1} holds for $f^k$. Hence, if $\tilde{x}\in I$ is an attracting periodic point of $f^k$ then in its immediate basin of attraction $B(\tilde{x})$ there is either an endpoint of $I$ or a critical point $\tilde{c}$ of $f^k$. Since any point that is attracted to the orbit of $\tilde{x}$ under $f^k$ is also attracted to the orbit of $\tilde{x}$ under $f$ the immediate basin of attraction of $\tilde{x}$ under $f$ contains $B(\tilde{x})$ hence either an endpoint of $I$ or $\tilde{c}$. As the critical set of $f^k$ is given by $$\mathcal{C}(f^k)=\{x\in I: \ \exists  \ \  0\leq i \leq k-1 \ \text{where} \ f^i(x)\in \mathcal{C}(f)\}$$ then $\tilde{c}$ is the preimage of some critical point of $f$ and the first and second statement of theorem \ref{theorem 1} follows from the fact that $f$ has an attracting periodic orbit containing $\tilde{x}$ if and only if the same is true of $f^k$. Since this is true for general periodic points the result follows.
$\Box$

\begin{definition} Let $f : I \rightarrow I$ be a $C^1$ map with nonempty critical set $\mathcal{C}(f)$. Let $$a=\inf_{n\geq 0}\{f^n(c) : \ c\in \mathcal{C}(f)\}, b=\sup_{n\geq0}\{f^n(c) : c\in\mathcal{C}(f)\}.$$
We call the interval $[a,b]$ the \textit{critical interval} of the function $f$. That is, the smallest closed interval that contains the forward orbit of all critical points of $f$.
\end{definition}

By a simple argument it follows that if $f$ has a critical interval $\tilde{I}$ then $f(\tilde{I})\subseteq\tilde{I}$. A more complicated result is the following.

\begin{lemma}
Let $f : I \rightarrow I$ be $C^1$ with a nonempty critical set. Then the endpoints of the critical interval of $f$ attract some critical point of $f$. Specifically, these endpoints are either attracting periodic points of period 1 or 2 or else lie on the orbit of a critical point of $f$.
\label{lemma 3.3}
\end{lemma}

\noindent\textit{Proof.}
For simplicity let $I=[0,1]$ which can always be achieved by some affine change of coordinates. Let $I=L\cup C \cup R$ where $C=[c_l,c_r]$ is the smallest interval containing $\mathcal{C}(f)$, $L=[0,c_l]$, and $R=[c_r,1]$. Also $c_{max}$ and $c_{min}$ are the critical points of $f$ with largest and smallest function values respectively. For the critical interval $[a,b]$ suppose in the following cases that $b$ is not on the orbit of any critical point. 

\textit{Case 1}: Let $f$ be increasing on $L$ and $R$. If $f^2(c_{max})\leq f(c_{max})$ then it follows that $f([0,f(c_{max})])\subseteq [0,f(c_{max})]$ so either $b=f(c_{max})$ or $b=c_r$, which violates the supposition. As it follows then that $f^2(c_{max})> f(c_{max})$ then $c_r<f(c_{max})$ implying that $f$ is strictly increasing on $[f^2(c_{max}),1]$. Hence, by monotonicity any orbit containing a point in $[f^2(c_{max}),1]$ is attracted to a fixed point of this interval. The least of these fixed points $p$ must be $b$ since $c_{max}$ is attracted to it and as $f([0,p])\subseteq [0,p]$. 

\textit{Case 2}: Let $f$ be decreasing on $L$ and $R$. Consider $f^2$ where we define $L^2$, $C^2$, $R^2$, $c^2_{max}$ analogous to $L,C,R$ and $c_{max}$ for $f^2$. Note that $f^2(c^2_{max})=f(c_{max})<b$. Also if $f^2$ is decreasing on $L^2$ this implies that $f(L)\subseteq [c_l,c_r]$ and $f^2$ decreasing on $R^2$ implies that $f(R)\subseteq [c_l,c_r]$. If either of these is the case then either $f(c_{max})$ or $f^2(c_{min})$ is equal to $b$. As neither of these is possible then $f^2$ is increasing on both $L^2$ and $R^2$ and the analysis reduces to that of Case 1. Therefore, $c^2_{max}$ is attracted to $b$ which implies the critical point $c_{max}$ of $f$ is attracted to a two cycle of $f$ containing $b$. 

\textit{Case 3}: Let $f$ be decreasing on $L$ and increasing on $R$. If $f^2(c_{max})\leq f(c_{max})$ and $f^2(c_{min})\leq f(c_{max})$ then $f([f(c_{min}),f(c_{max})])\subseteq [f(c_{min}),f(c_{max})]$ implying $b$ lies on the orbit of $c_{max}$ or $c_r$, which violates the supposition. If $f^2(c_{max})> f(c_{max})$ then as in Case 1 either $c_{min}$ or $c_{max}$ is attracted to a fixed point in the interval $[f^2(c_{max},1)]$ which must be $b$. If $f^2(c_{min})>f(c_{max})$ but $f^2(c_{max})\leq f(c_{max})$ then $f^2(c_{min})=b$.

\textit{Case 4}: Let $f$ be increasing on $L$ and decreasing on $R$. This however implies that $b=f(c_{max})$.

Repeating this argument with appropriate modifications implies the same is true for the endpoint $a$. 
$\Box$

\begin{corollary}\label{corollary 1}
An $S^k$-multimodal map $f$ can have at most $|\mathcal{C}(f)| + 2$
attracting periodic orbits. If no critical point of $f$ is attracted to a periodic orbit then all periodic points in the critical interval are hyperbolic repelling.
\end{corollary}

\noindent\textit{Proof.}
The first statement is immediate from theorem \ref{theorem 1}. To prove the second note that lemma \ref{lemma 3.3} implies $f^k$ restricted to its critical interval $\tilde{I}$ has the property that each of its attracting periodic orbits attracts a critical point. Hence, if no critical point of $f^k$ is attracted to a periodic orbit, $f^k$ can have no attracting periodic points in its critical interval, in particular no hyperbolic attractors as well as no neutral periodic points in the interior of $\tilde{I}$. If an endpoint of $\tilde{I}$ is a neutral periodic point not on the orbit of a critical point it attracts an open set of points in the critical interval and the corollary follows for $f^k$. Hence, all periodic orbits of $f$ in $\tilde{I}$ are hyperbolic repelling and as the critical interval of $f$ is contained in $\tilde{I}$ the proof follows.
$\Box$

\section{Measure Theoretic Properties} 

In this section we prove theorems $\ref{theorem 2}$, $\ref{theorem 3}$, and $\ref{theorem 4}$ from section 1, that is those results that are more measure theoretic in nature. We also give some related corollaries that will be used in section 6 in our discussion of neuronal models. We are now give the proof of theorem $\ref{theorem 2}$.

\noindent\textit{Proof.}
For any $\tilde{c}\in \mathcal{C}(f^k)$ only one of $\tilde{c},f(\tilde{c}),f^2(\tilde{c}),\dots,f^{k-1}(\tilde{c})$ can be a critical point of $f$ since if $f$ eventually maps any point of $\mathcal{C}(f)$ back to this set then condition (\ref{eq:no2}) does not hold. Hence, for $\tilde{c}\in\mathcal{C}(f^k)$ there exists a unique $m<k$ such that $f^m(\tilde{c})=c$ where $c\in\mathcal{C}(f)$. For the sake of convenience let 
\begin{align*}
D_{n,k}(\tilde{c})&=|((f^k)^n)^{\prime}(f^k(\tilde{c}))|,\\
C_{k,m}(\tilde{c})&=\prod_{i=1}^{k-m-1}|f^{\prime}(f^i(\tilde{c}))|, \ m<k.
\label{eq:no6}
\end{align*} 
By repeated use of the chain rule it follows that 

\begin{equation*} 
C_{k,m}(\tilde{c})D_{n,k}(\tilde{c})=\prod_{i=1}^{nk-m-1}|f^{\prime}(f^i(c))|.
\end{equation*} The important observation here is $C_{k,m}(\tilde{c})\neq 0$ and does not depend on $n$ and also $D_{nk-m-1}(c)=C_{k,m}(\tilde{c})D_{n,k}(\tilde{c})$ for $n\geq 1.$ Since all the quantities involved are positive
\begin{equation}
C_{k,m}^{-1/(2\ell-1)}(\tilde{c})\sum_n D_{n,k}^{-1/(2\ell-1)}(\tilde{c})\leq \sum_n D_{n}^{-1/(2\ell-1)}(c) < \infty. 
\label{eq:no7}
\end{equation}

Let $\tilde{x}$ be either a hyperbolic attracting or neutral fixed point of $f^k$ in its critical interval. Theorem \ref{theorem 1} together with lemma \ref{lemma 3.3} imply that some critical point $\tilde{c}\in \mathcal{C}(f^k)$ is attracted to $\tilde{x}$. In either case there is an $N$ and a closed interval $J\ni\tilde{x}$ on which $|(f^k)^\prime(x)|\leq 1$ containing $f^{kn}(\tilde{c})$ for $n\geq N$. Hence, there is a finite $C\geq 0$ such that for $n\geq N$, $D_{n,k}(\tilde{c})\leq C$ implying
$$\sum_n D_{n,k}^{-1/(2\ell-1)}(\tilde{c})\geq\sum_{n\geq N}C^{-1}=\infty.$$
As this violates (\ref{eq:no7}) it follows that every fixed point in the critical interval of $f^k$ is hyperbolic repelling. A slight modification of this argument implies the same is true for all periodic orbits in the critical interval of $f^k$ as well. 

Note that the critical interval of $f$ is not a single point since a critical point of $f$ cannot be a fixed point of the function if (\ref{eq:no2}) holds. Hence, as the critical interval of $f^k$ contains that of $f$, if $\tilde{f}$ denotes $f$ restricted to its critical interval then the main result of \cite{8} implies the result of the theorem for $\tilde{f}$. This result can then be trivially extended to $f$.
$\Box$\\

For the proof of theorem \ref{theorem 3} we require the following lemma.

\begin{lemma}
Let $f$ be $C^1$ with a finite critical set. If for some $k\in\mathbb{N}$ there is an $f^k$-invariant probability measure $\mu$ absolutely continuous with respect to Lebesgue measure (acip) then $f$ admits an invariant measure also absolutely continuous with respect to Lebesgue measure.
\label{lemma 4.2}
\end{lemma}

\noindent\textit{Proof.}
Suppose $f,k,$ and $\mu$ are as above. Consider the measure $\upsilon$ given by 

$$\upsilon(A)=\frac{1}{k}\sum_{i=0}^{k-1}\mu(f^{-i}(A)).$$
To see that $\upsilon$ is $f$-invariant note 

$$\upsilon(f^{-1}(A))=\frac{1}{k}\sum_{i=1}^{k-1}\mu(f^{-i}(A))+\frac{\mu(f^{-k}(A))}{k}=
\frac{1}{k}\sum_{i=1}^{k-1}\mu(f^{-i}(A))+\frac{\mu(A)}{k}=\upsilon(A).$$ 
For absolute continuity of the measure note that a $C^1$ function with a finite critical set is non-singular. Therefore, absolute continuity of the measure $\upsilon$ follows from that of $\mu$.
$\Box$\\ 

We now give the proof of theorem $\ref{theorem 3}$

\noindent\textit{Proof.}
Assuming condition (\ref{eq:no3}) if $\tilde{c}\in \mathcal{C}(f^k)$ then, as in the proof of theorem \ref{theorem 2}, there is a unique $m<k$ such that $f^m(\tilde{c})=c$ where $c\in \mathcal{C}(f)$. Suppose $c$ has order $\ell$ under $f$. Consider 
\begin{equation}
\lim_{x\rightarrow \tilde{c}}\frac{|(f^k)^{\prime}(x)|}{|x-\tilde{c}|^{\ell-1}}=\prod_{i=0,i\neq m}^{k-1}|f^{\prime}(f^i(\tilde{c}))|\lim_{x\rightarrow \tilde{c}}\frac{|f^{\prime}(f^m(x))|}{|x-\tilde{c}|^{\ell-1}}.
\label{eq:no4}
\end{equation}
Let $A=\prod_{i=0,i\neq m}^{k-1}|f^{\prime}(f^i(\tilde{c}))|$ which is strictly positive and note that 
$$\lim_{x\rightarrow \tilde{c}}\frac{|f^{\prime}(f^m(x))|}{|f^m(x)-f^m(\tilde{c})|^{\ell-1}}=\lim_{x\rightarrow f^m(\tilde{c})}\frac{|f^{\prime}(x)|}{|x-f^m(\tilde{c})|^{\ell-1}}>0$$ where the inequality follows from definition 2.3(ii). Setting this limit to $B$, the right side of equation ($\ref{eq:no4}$) is
$$A\lim_{x\rightarrow \tilde{c}}\frac{|f^{\prime}(f^m(x))|}{|f^m(x)-f^m(\tilde{c})|^{\ell-1}}\frac{|f^m(x)-f^m(\tilde{c})|^{\ell-1}}{|x-\tilde{c}|^{\ell-1}}=AB\lim_{x\rightarrow \tilde{c}}\Big(\frac{|f^m(x)-f^m(\tilde{c})|}{|x-\tilde{c}|}\Big)^{\ell-1}.$$ 
An application of L'Hospital's rule implies that this limit is strictly positive. That is, the collection of orders of $\mathcal{C}(f)$ is the same as those for $\mathcal{C}(f^k)$ implying their maximum $\ell_{max}$ is equal. 

By replacing $-1/(2\ell-1)$ with $-1/\ell_{max}$ in the derivation of (\ref{eq:no7}) we obtain $$\sum_n D_{n,k}^{-1/(\ell_{max})}(c)< \infty \ \text{for all}  \ c\in \mathcal{C}(f^k).$$ 
This implies $f^k$ has an acip via the main result in \cite{4} and lemma \ref{lemma 4.2} implies the same for $f$.
$\Box$\\
 
We now proceed to the proof of theorem $\ref{theorem 4}$.

\noindent\textit{Proof.} 
As each point in $\mathcal{C}(f_{\beta}^k)$ is the preimage of some critical point $c$ of $f_{\beta}$ then inf$_{i>0}d((f_{\beta}^k)^i(c),\mathcal{C}(f_{\beta}^k))>0$ for all $c\in \mathcal{C}(f_{\beta}^k)$ as this would otherwise violate condition (iv) on $f$. Second, for $c\in \mathcal{C}(f_\beta)$ note that $(f_{\beta})^{\prime\prime}(c)\neq0$ implies $c$ is nonflat with order $\ell=2$. As inf$_{i>0}d((f_{\beta}^k)^i(c),\mathcal{C}(f_{\beta}^k))>0$ for all $c\in \mathcal{C}(f_{\beta}^k)$ then the same argument used in the proof of theorem \ref{theorem 3} implies that the critical points of $f^k_\beta$ are also nonflat of order $\ell=2$. Therefore, $(f^k_\beta)^{\prime\prime}(c)\neq 0$ for all $c\in \mathcal{C}(f^k_\beta)$. Also note that property (iii) of theorem \ref{theorem 4} means $f_{\beta}(x)$ has no attracting or neutral periodic orbit. But this is true if and only if $f_{\beta}^k(x)$ has none itself. 

From the assumption that $f^k_{\beta}$ has a negative Schwarzian derivative it follows from \cite{9} that $f^k_{\beta}\in\mathcal{M}$. In particular, $f^k$ admits an invariant absolutely continuous measure $\mu$. And an application then of lemma \ref{lemma 4.2} implies that $f_{\beta}$ also has an acip. Moreover, if $f^k_{\beta}$ satisfies (\textit{PT}) then on a positive Lebesgue measure set of parameters the family of functions $f^k_a$ has an acip and lemma \ref{lemma 4.2} can again be used to show the same for $f_a$.
$\Box$\\
 
As we will be concerned specifically with unimodal maps in section 6 we give the following corollary.

\begin{corollary}\label{corollary 2}
Let $f:I \rightarrow I$ be $C^3$ and $S^k$-unimodal with critical point $c$ of order $\ell=2$ to the left of which $f$ is increasing and to the right of which is decreasing. If the orbit of $c$ contains a repelling periodic orbit and $f$ has no fixed points on the boundary or outside of its critical interval then $f^k\in\mathcal{M}$.
\end{corollary}

\noindent\textit{Proof.}
Assuming these conditions then $f(c)>c$ since $c$ is otherwise attracted to the rightmost hyperbolic attracting or neutral fixed point of $f$. Also if $f^2(c)\geq c$ the forward orbit of $c$ is contained in the interval $[c,f(c)]$ so $c$ is attracted to a fixed point or a periodic cycle neither of which can be repelling. Similarly, if $f^3(c)<f^2(c)$ then $c$ is also attracted to a nonrepelling fixed point, this implies that the critical interval $\tilde{I}$ of $f$ is $[f^2(c),f(c)]$.   

Note that if $f$ has no fixed points on the boundary or outside $\tilde{I}$ then there is some $n\in\mathbb{N}$ such that for all $x\in I\setminus\tilde{I}$, $f_{\beta}^n(x)\in\tilde{I}$ since $\tilde{I}$ is forward invariant and all orbits fall on this set after some bounded number of iterations. Therefore, $f$ has no periodic points outside $\tilde{I}$ and the proof of theorem \ref{theorem 2} implies that on $I$, $f$ satisfies condition (iii) of theorem \ref{theorem 4} as $c$ is mapped to an unstable periodic orbit. Condition (iv) follows for the same reason, (ii) follows from the proof of theorem \ref{theorem 4} and condition (i) is assumed to hold. 
$\Box$

\section{Characterizing Functions with an Eventual Negative Schwarzian Derivative}

As not every $C^2$ function has the characteristics given in theorem \ref{theorem 1} it is not the case that every function will be either a function with a negative Schwarzian derivative or have an eventual negative Schwarzian derivative. That is, it is possible for a function to have a Schwarzian derivative that is mixed for all of its iterates. If however a function does have an eventual negative Schwarzian derivative it would be useful to have a way of identifying this. Specifically, we would like to have sufficient conditions under which a function has this property.  

The simplest case to consider is the one in which $|f^\prime|^{-1/2}$ is convex but not strictly convex. To do so we mention the following. 

\begin{definition} 
Let $g:J\rightarrow K$ be continuous and monotone where $U\subset V\subset J$ are open bounded intervals in $K$. If $V\setminus U$ consists of the intervals $L$ and $R$ the \textit{cross ratio} of the intervals $U$ and $V$ is given by $$CR(U,V)=\frac{|U||V|}{|L||R|}.$$ The function $g$ is said to expand cross ratios on $J$ if for any intervals $U\subset V$ in $J$, $CR(U,V)<CR(g(U),g(V)).$  
\end{definition}

Suppose $f$ is $C^2$ on the open, bounded interval $J$ containing no critical points of $f$. It is known that $f$ increases cross ratios on $J$ if and only if it has a negative Schwarzian derivative on this set (see \cite{12}). Furthermore, if a function is a M\"{o}bius transform, that is a function of the form $g(x)=(ax+b)/(cx+d)$ where $ad-bc\neq 0$, then this function preserves cross ratios i.e. $CR(U,V)=CR(g(U),g(V))$ (see \cite{14}).  

\begin{proposition}\label{proposition 1}
Let $f:I\rightarrow I$ be $C^2$. Suppose $M\subset I$ is a finite union of closed intervals on which $f$ is a M\"{o}bius transform and off of which $f$ has a negative Schwarzian derivative. Furthermore, assume there is a $k\geq 2$ such that for every $x\in M$, $\{f^{i}(x): 1\leq i < k\}\cap I\setminus M \neq \emptyset$. Then $f$ has an eventual negative Schwarzian derivative of order less than or equal to $k$.
\end{proposition}

\noindent\textit{Proof.}
If $B$ are the boundary points of $M$ let $B^k=\{x\in I: f^{i}(x)\in B\ \ \text{for some} \ 0\leq i< k\}$. Let $J\subset I$ be an open interval containing no points in $\mathcal{C}(f^k)$ or $B^k$. Then there is a first $0\leq \ell< k$ such that $f^\ell(J)\subset I\setminus M$. As the composition of a M\"{o}bius transform with itself is again such a function then $f^\ell(x)$ restricted to $J$ does not change cross ratios. However, $f^\ell(J)\cap \mathcal{C}(f)=\emptyset$ and $f^\ell(J)\subseteq I\setminus M$ on which $f$ has a negative Schwarzian derivative implying $f^{\ell+1}$ increases cross ratios on $J$. Since cross ratios on $J$ are either maintained or increased by further iteration of $f$ it follows that $f^{k}$ also increases cross ratios on $J$ or $|(f^{k})^\prime|^{-1/2}$ is strictly convex on $J$. Note that as $|(f^{k})^\prime|^{-1/2}$ is $C^1$ on $J$ it then has a strictly increasing derivative on this interval. Since $B^k$ is a finite set it follows that on any open interval of $I$, not containing critical points of $f^k$, $|(f^{k})^\prime|^{-1/2}$ is also strictly convex or $f^k$ has a negative Schwarzian derivative on $I$. 
$\Box$

It follows directly from this proposition that the $C^2$ family of unimodal functions 
\begin{equation} 
 f_a(x)=\begin{cases}
  (x-1/2)+a, & x\in [0,1/2]\\
  -4(2a+1)(x-1/2)^3+(x-1/2)+a, & x\in (1/2,1]\\
\end{cases}
\label{eq:no6}
\end{equation} has an eventual negative Schwarzian derivative for $a\in (1/2,7/8]$ (see Fig. 1). A $C^3$ example of a similar family is given by $g_a(x)=(ax-1/4)+(ax-1/4)^4+4a/11$ for $x\in[0,1]$ and $a\in[1,11/8]$. The reason for (\ref{eq:no6}) is that it serves as a simple example of a family of functions for which many results, which refer only to $C^3$ functions, are not applicable. For example, as $f_a(x)$ is only $C^2$ it is not directly possible to use the main results in \cite{5} to show this family is conjugate to a function with a negative Schwarzian derivative. 

It is also worthwhile to recall, as mentioned in the introduction, the motivation as well as the inspiration for the study of this new class of functions with eventual negative Schwarzian derivatives comes from the analysis of one-dimensional maps which appear in some models in neuroscience. These maps, given in \cite{10,11} in particular, have a part that is linear (or almost linear). Therefore, our example above also includes this feature although, as implied by proposition \ref{proposition 1}, this is not a necessary condition to have an eventual negative Schwarzian derivative.

\begin{figure}
\begin{center}
\begin{overpic}[scale=0.6]{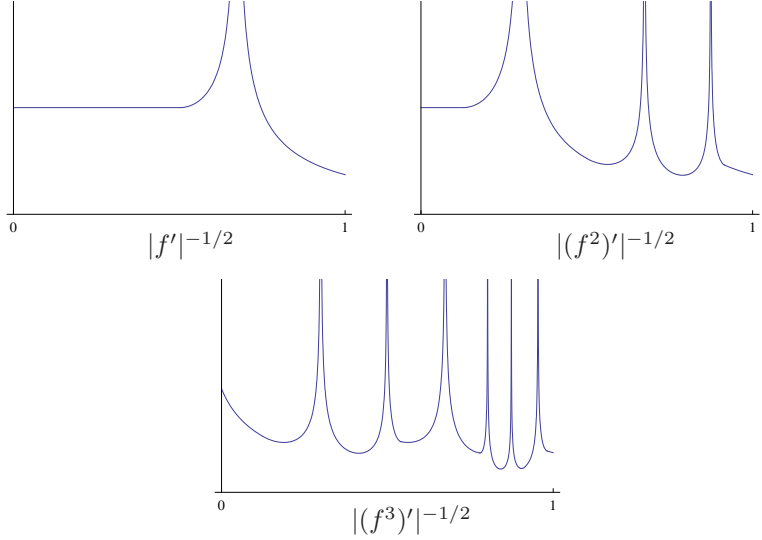}
\put(20,33){$|f^\prime|^{-1/2}$}
\put(72,33){$|(f^2)^\prime|^{-1/2}$}
\put(46,-2){$|(f^3)^\prime|^{-1/2}$}
\end{overpic}
\end{center}
\caption{$f=f_{7/8}$ in equation (\ref{eq:no6})}.
\end{figure} 
 
From (\ref{eq:no8}) it follows that if $f$ is $C^3$ then
\begin{equation}
S(f^k)(x)=\sum_{i=0}^{k-1} \Big((f^{i})^{\prime}(x)\Big)^2 \Big(S(f)(f^i(x))\Big).
\label{eq:no9}
\end{equation}

Note that if in the previous proposition the function was assumed $C^3$ with $S(f)(x)<0$ on $I\setminus M$, this equation would have immediately implied the result. However, this equation suggests a method for identifying $C^3$ functions which have an eventual negative Schwarzian derivative. 

If $I=\bigcup_{i=1}^n I_i$, where the $I_i$ are nonintersecting intervals, let $A=(a_{ij})$ be the \textit{transition matrix} of this \textit{partition} with respect to $f$. That is, $a_{ij}=1$ if there is an $x\in I_i$ such that $f(x)\in I_j$ and $a_{ij}=0$ otherwise. Let a sequence $\bar{x}=(x_0x_1x_2\dots x_{k-1})$ of length $k$ be \textit{admissible} with respect to this partition if each $x_i\in\{1,2,\dots,n\}$ and $x_j$ can follow $x_i$ in this sequence if and only if $a_{ij}=1$. Let $T(i)=\sup_{x\in I_i}\{S(f)(x)\}$, $m(i)=\inf_{x\in I_i}\{|f^{\prime}(x)|^2\}$, $M(i)=\sup_{x\in I_i}\{|f^{\prime}(x)|^2\}$ and for some admissible $\bar{x}$ of length $k$ and $0<j<k$ define 
$$R_j(\bar{x})=\begin{cases}
  \big[\prod_{i=0}^{j-1}m(x_i)\big]T(x_j), & T(x_j)\leq 0\\
  \big[\prod_{i=0}^{j-1}M(x_i)\big]T(x_j), & \text{otherwise}\\
\end{cases}$$ 
letting $R_0(\bar{x})=T(x_0)$. From equation (\ref{eq:no9}) we have the following proposition.
    
\begin{proposition}
Let $f:I\rightarrow I$ be $C^3$ and $I=\bigcup_{i=1}^n I_i$ a partition of $I$. If there is a $k\geq 1$ such that for every admissible sequence $\bar{x}$ of length $k$, 
$$\sum_{i=0}^{k-1}R_i(\bar{x})<0$$ then $f$ has an eventual negative Schwarzian derivative of order less than or equal to $k$.  
\end{proposition}

Using this proposition it can be shown that the one parameter family of functions 
\begin{equation}
g_a(x)=1-a\tan\Big(\frac{\pi}{4}x^2\Big), \ x\in [-1,1], 
\label{eq:no10}
\end{equation} which mimics the logistic function, has an eventual negative Schwarzian derivative of order $k=2$ in a parameter neighborhood of $a=1.7$. This can be done by using the partition with endpoints $\{-1,-0.95,-0.7,-0.47,-0.18,0.18,0.47,0.7,0.95,1\}$ for instance.  

\section{Application to a Neuronal Model}
The motivation for considering functions having some iterate with a negative Schwarzian derivative comes from a model for the electrical activity in neural cells specifically in behavior described as bursting. This model given first in \cite{10} and later in \cite{11} is a reduction of a system of three nonlinear differential equations to a 1-$d$ map. 

The model is initially given by the following fast-slow system of three differential equations which describe the dynamics of the membrane potential $\upsilon$ and two gating variables $\eta$ and $\omega$ of a neuronal cell:
\begin{align} 
\epsilon\dot{\upsilon}&=f(\upsilon, \eta, \omega;\delta)\\
\dot{\eta}&=g(\upsilon,\eta)\\
\dot{\omega}&=\beta h(\eta,\omega)
\end{align}
Here the parameter $\delta$ can be viewed as a control parameter of the full system where $\delta\in[\delta_{min},\delta_{max}]$. Also the time constant $\beta$ represents the slowest time scale in the dynamics of (10)-(12) so in the limit $\beta\rightarrow 0^+$ the system uncouples into a fast subsystem (10), (11) and a slow subsystem (12). As is explained in \cite{11} the trajectory of the full system is drawn towards a surface foliated by periodic orbits of the fast subsystem where the evolution along this surface is determined by the dynamics of the slow subsystem. 

As the state of the fast subsystem depends on the value of the slow variable $\omega$, it is sufficient to know how $\omega$ changes after each oscillation of the fast subsystem. Knowing these changes is precisely the reduction of the system to a 1-$d$ map denoted by $F_\delta$. To achieve this a Poincar\'{e} section $P_\delta$ is placed transversal to the surface of periodic orbits and the map is defined by $F_\delta(\omega_{n})=\omega_{n+1}$ where $\omega_{n+1}$ is the $\omega$-coordinate of the next point on the flow to pass through $P_\delta$.

In \cite{10} $F_\delta$ is shown to have the following properties (see Fig. 2):

\newcounter{Lcount}
\begin{list}{\Roman{Lcount}}
{\usecounter{Lcount}
\setlength{\rightmargin}{\leftmargin}}
\item For fixed $\delta$, $F_\delta(x)$ is a piecewise $C^2$ map with two intervals of continuity, $I_1=I^-\cup I^0$ and $I_2=I^+$, between which there is a single discontinuity $d(\delta)$ continuous in $\delta$.
\item $F_\delta$ is unimodal on $0\leq x< d(\delta)=I_1$, with critical point $c(\delta)$ of order $\ell=2$.
\item $F_\delta$ is nearly linear on its left outer region $I^-$ with slope slightly less than 1 and in the limit $\epsilon\rightarrow 0^+$ is linear with slope $1$.
\item On the inner region $I^0$, $F_\delta$ is unimodal with slope tending toward $-\infty$ as $d(\delta)$ is approached from the left. 
\item In the second region $I^+$ the function is nearly constant and between 0 and $F_\delta(c(\delta))$.
\item $F_\delta$ has a unique fixed point $\alpha(\delta)$ on the interval $I^0$.
\item There is a $\delta_0\in (\delta_{min}, \delta_{max})$ where $\alpha(\delta_0)=c(\delta_0)$ and for $\delta\in (\delta_0,\delta_{max}]$, $c(\delta)<\alpha(\delta)$.
\item There is a $\delta_b\in (\delta_0,\delta_{max})$ such that $F_{\delta_b}(c(\delta_b))=d(\delta_b)$ and for $\delta\in [0,\delta_b)$, $F_\delta(c(\delta))<d(\delta)$.
\item There is a $\delta_n>\delta_0$ such that for $\delta\in(\delta_n,\delta_b)$, $F^\prime_\delta(\alpha(\delta))<-1$. 
  \end{list}

\begin{figure}
\begin{center}
\begin{overpic}[scale=.5]{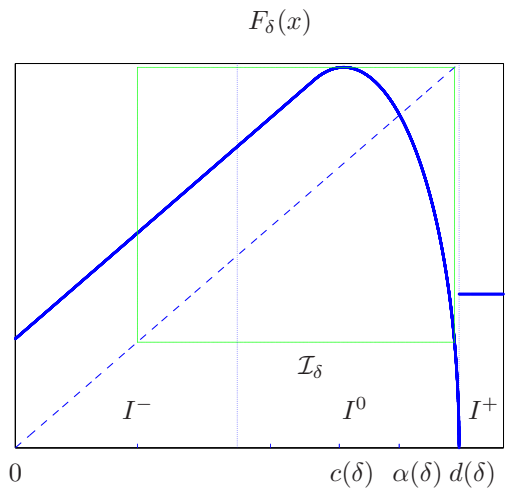}
\put(63,3){$c(\delta)$}
\put(73,3){$\alpha(\delta)$}
\put(58,20){$\mathcal{I}_{\delta}$}
\put(82,3){$d(\delta)$}
\put(12,3){0}
\put(50,75){$F_{\delta}(x)$}
\put(30,13){$I^-$}
\put(65,13){$I^0$}
\put(85,13){$I^+$}
\end{overpic}
\end{center}
\caption{First Return Map Near Bursting $\delta\approx\delta_b$}
\end{figure} 

\begin{remark}
As the map $F_\delta$ depends on the placement of $P_\delta$ then locally $F_\delta$ also corresponds to this placement. This variability in the placement of $P_\delta$ is one of the main obstacles in verifying whether $F_\delta$ has an eventual negative Schwarzian derivative. However, at a global level the nearly linear part of $F_\delta$ on $I^-$ and the fact that orbits leave this interval in a finite number of iterations suggests that if off $I^-$ the function has a negative Schwarzian derivative then for $\epsilon$ small enough $F_\delta$ will have an eventual negative Schwarzian derivative (see (\ref{eq:no6}) for the linear case). 
\end{remark}

Our aim in this section is to give sufficient conditions under which this one parameter family of functions has an acip for a positive Lebesgue measure set of its parameters. Also under what conditions these functions exhibit a mixing property with respect to these measures. These conditions will ultimately include the assumption that $F_\delta$ has an eventual negative Schwarzian derivative in order to illustrate the usefulness of this concept.

To simplify notation denote $F^n_\delta(c(\delta))=c^n(\delta)$, $F_\delta$ restricted to its critical interval by $\tilde{F_\delta}$, and let $\mathcal{O}=\{(\delta,x): \delta\in(\delta_0,\delta_b), \ x\in(0,d(\delta))\}$. In what follows we will often consider $F_\delta(x)=F(\delta,x)$ to be a function of two variables on $\mathcal{O}$. It will be the open set $\mathcal{O}$ on which we will focus our attention. The reason being is that for parameter values $\delta\leq\delta_0$ the map $F_\delta$ has a global attracting fixed point and for parameters larger than $\delta_b$ the map $F_\delta$ is not continuous so the previous theory does not apply. Specifically, we make the observation that so long as $\delta\in(\delta_0,\delta_b)$ then $c^1(\delta)<d(\delta)$ implying the critical interval $\mathcal{I}_{\delta}$ of $F_\delta$ is either\\
(i) $[c(\delta),c^1(\delta)]$ if $c(\delta)\leq c^2(\delta)$ in which case $\tilde{F_\delta}$ is a diffeomorphism or\\
(ii) $[c^2(\delta),c^1(\delta)]$ if $c(\delta)>c^2(\delta)$ and  $\tilde{F_\delta}$ is unimodal. 

The following lemma is meant to establish basic continuity properties of the family of functions $F_\delta$ under variation of parameters.

\begin{lemma}\label{lemma 6.1} Suppose $F_\delta(x)=F(\delta,x)$ is $C^2$ on $\mathcal{O}$ and there exist a closed interval $J$ with endpoints $\delta_1,\delta_2\in (\delta_0,\delta_b)$ such that for all $\delta\in J$, $c^3(\delta)<\alpha(\delta)$. If there is also an $m\in\mathbb{N}$ where $c^m(\delta_1)= c^1(\delta_1)$ and $c^m(\delta_2)= \alpha(\delta_2)$ then on an infinite set of parameters $\Delta\subset J$ the orbit of $c(\delta)$ contains $\alpha(\delta)$.  
\end{lemma}

\noindent\textit{Proof.}
If $\delta\in J$ then $c(\delta)>c^2(\delta)$ implying $\mathcal{I}_\delta=[c^2(\delta),c^1(\delta)]$ on which $F_\delta$ is unimodal. Also every point in $(\alpha(\delta), c^1(\delta))$ has exactly two preimages; one preimage in $l(\delta)=(c^2(\delta),c(\delta))$ and the other in $r(\delta)=(c(\delta),\alpha(\delta))$. As every point in $l(\delta)$ and $r(\delta)$ has a unique preimage in $s(\delta)=(\alpha(\delta), c^1(\delta))$ it is possible to specify a preimage of $\alpha(\delta)$ by some finite sequence made up of $l,r$, and $s$ which stand for whether the fixed point is reached by tracing its forward orbit through these sets in this manner. We say a finite sequence of $l,r,$ and $s$ is admissible if and only if an $s$ separates every $l$ or $r$, $s$ does not follow itself, and the sequence ends with an $l$. This corresponds to the structure above and uniquely defines a preimage of $\alpha(\delta)$. Note we need not assign a symbol to $c(\delta),$ $c^1(\delta),$ $c^2(\delta)\in\mathcal{I}_\delta$ since no preimage of $\alpha(\delta)$ with an admissible sequence has an orbit containing these points for $\delta\in J$.

Under the assumption that $F_\delta(x)=F(\delta,x)$ is $C^2$ the Implicit Function Theorem guarantees that $c(\delta)$ is a $C^1$ function of $\delta$ on $J$ as the critical point of $F_\delta$ is to leading order quadratic. By similar calculations it follows that the same is true for $\alpha(\delta)$, and $c^m(\delta)$ for all $m\geq 1$. Also important is that this is also true for any preimage of $\alpha(\delta)$ having an admissible sequence. This follows as well from the Implicit Function Theorem using the fact that the orbit of these particular preimages cannot contain the critical point.

As there are infinitely many sequences which can end in either $l,r,s$ then infinitely many admissible preimages of the fixed point $\alpha(\delta)$ are in each of $l(\delta)$, $r(\delta)$, and $s(\delta)$ and each vary continuously in $\delta$. As the first letter in the sequence of a preimage determines whether it is in $l(\delta),$ $r(\delta),$ or $s(\delta)$ it must stay in that interval for all $\delta\in J$.

Continuity, specifically of the admissible preimages in $s(\delta)$, and the assumption that $c^m(\delta_1)= c^1(\delta_1)$ and $c^m(\delta_2)=\alpha(\delta_2)$ then imply the existence of the infinite set $\Delta\subseteq[\delta_1,\delta_2]$ on which the orbit of the critical point $c(\delta)$ contains the fixed point $\alpha(\delta)$. 
$\Box$

\begin{theorem}\label{theorem 6}
Suppose $F_\delta(x)=F(\delta,x)$ is $C^2$ on $\mathcal{O}$ and assume that there are $\delta_1,\delta_2\in(\delta_0,\delta_b)$ and $m\geq 3$ such that\\
(i) $c(\delta_1)\leq c^m(\delta_1)$,\\
(ii) $c^m(\delta_2) \leq c(\delta_2)$ and $c^{3}(\delta_2),$ $c^{m+1}(\delta_2)< \alpha(\delta_2)$.\\
Then there is an infinite set $\Delta\subset[\delta_1,\delta_2]$ on which the orbit of $c(\delta)$ contains $\alpha(\delta)$. Furthermore, if for some $\delta\in\Delta$\\
(iii) $\tilde{F_\delta}$ has an eventual negative Schwarzian derivative of order $k$ and\\
(iv) $\delta>\delta_n$ then the following is true:\\
\indent(A) $\tilde{F}_\delta$ has an acip.\\
\indent(B) If $\tilde{F_\delta}$ is $C^3$ in $x$ then some iterate of $\tilde{F_\delta}$ is mixing with exponential decay of correlations. 
\end{theorem}

\noindent\textit{Proof.}
As in the proof of lemma \ref{lemma 6.1} the assumption that $F_\delta$ is $C^2$ implies $\alpha(\delta)$, $c(\delta)$, and $c^i(\delta)$ for $i\geq 1$ all vary continuously in the parameter. Without loss in generality if $\delta_1<\delta_2$ since both $c(\delta_1)\leq c^{m}(\delta_1)$ and $c^m(\delta_2)\leq c(\delta_2)$ then there is a largest $\delta_*\in [\delta_1,\delta_2]$ such that $c(\delta_*)=c^m(\delta_*)$. As this implies $c^1(\delta_*)=c^{m+1}(\delta_*)$ then $\delta_*\neq \delta_2$. Note that if at some $\beta\in[\delta_*,\delta_2]$, $c^3(\beta)\geq \alpha(\beta)$ the assumption $c^3(\delta_2)<\alpha(\delta_2)$ implies there is a $\gamma\in[\beta,\delta_2]$ at which $c^3(\gamma)=\alpha(\gamma)$ in turn implying $c^m(\gamma)=\alpha(\gamma)>c(\gamma)$ contradicting the maximality of $\delta_*$. Lemma \ref{lemma 6.1} therefore guarantees the existence of the set $\Delta\subset[\delta_*,\delta_2]$ as $c^{m+1}(\delta_*)=c^1(\delta_*)$ and $c^{m+1}(\delta_2)<\alpha(\delta_2)$.

If $\delta\in\Delta$ and condition \textit{(iv)} holds then there is $C>0$ such that $D_n(c(\delta))=C|\tilde{F}_\delta^{\prime}(\alpha(\delta))|^n$ for $n$ large enough where $|\tilde{F}_\delta^{\prime}(\alpha(\delta))|>1$. Hence, inequality (\ref{eq:no3}) of theorem \ref{theorem 3} holds. This together with condition \textit{(iii)} implies \textit{(A)} via theorem \ref{theorem 3}.

For \textit{(B)} if $\tilde{F_\delta}$ is $C^3$ in $x$ then from the calculation of $D_n(c(\delta))$ above theorem \ref{theorem 2} implies some iterate $\tilde{F_\delta}$ is mixing with exponential decay of correlations. 
$\Box$

\begin{corollary} 
For $\delta\in(\delta_0,\delta_b)$ let $h_{\delta}(x)=(b_{\delta}-a_{\delta})x+a_{\delta}$ where $[a_{\delta},b_{\delta}]$ is the critical interval of $F_\delta$. Also denote by $H_\delta:[0,1]\rightarrow[0,1]$ the family of functions $\tilde{F}_\delta$ conjugated by $h_\delta$. If $\tilde{F}_\delta$ satisfies conditions \textit{(i)-(iv)} of theorem \ref{theorem 6}, is $C^3$, and for some $\delta\in\Delta$, $H_\delta^k$ has property (\textit{PT}), then on a positive set of parameters $F_\delta$ has an acip.
\end{corollary}

\noindent\textit{Proof.}
Note that the one parameter family of functions $H_\delta=h^{-1}_{\delta}\circ F_{\delta}\circ h_{\delta}$ when restricted to the interval $[0,1]$ is a linearly scaled version of the family $F_\delta$ restricted to their critical intervals. It follows by a simple calculation that $S(H_\delta^k)(x)=S(F_\delta^k)(h_\delta(x))$ or the property of having an eventual negative Schwarzian derivative is preserved under this change of coordinates. Furthermore, $H^k_\delta$ has a finite number of nonflat critical points and as the orbit of $c(\delta)$ contains $\alpha(\delta)$ corollary \ref{corollary 2} implies that $H^k_\delta\in\mathcal{M}$. It follows directly from \cite{9} that if $H^k_\delta$ satisfies property    \textit{(PT)} then on a positive set of parameters $H_\delta$ has an acip which implies the same for $F_\delta$.
$\Box$\\

Numerically, as $\delta\rightarrow\delta_b$ the number of iterates the orbit of a point can stay on $\mathcal{I}_\delta \bigcap I^-$ increases from 1 near $\delta_0$ to around 5 near $\delta_b$ (see Fig.3 in \cite{11}). From this point of view conditions \textit{(i)} and \textit{(ii)} of the previous theorem are a natural way to ensure this family of maps have a parameter set $\Delta$ mentioned above. More than this, however, it allows for some latitude in selecting how close $\Delta$ is to $\delta_b$ which is useful since for parameters near $\delta_b$, the map $F_\delta$ is more likely to have a repelling fixed point. Again, because of the variability in the placement of $P_\delta$ it would take some work to verify whether or not for a particular placement condition \textit{(iii)} holds.

It should be noted that functions with the property that their critical points are mapped to repelling periodic orbits, often called \textit{Misiurewicz maps}, are very rare both in a topological and metric sense \cite{13}. From this the previous theorem seems to imply that mixing is a rare event for $F_\delta$. However, the conditions of theorem \ref{theorem 2} are generically satisfied on a larger set of parameters than just $\Delta$. The problem which is generally encountered is in showing for a particular parameter value or set of values that the derivatives along the orbits of the critical points have the requisite growth. However, one reason to expect that (B) holds for a much larger set of parameters than just $\Delta$ is the very large negative slope of $F_\delta$ near $d(\delta)$, for $\delta\approx \delta_b$.

\section{Conclusion}
The goal of this paper was first and foremost to study a more general class of functions that behave in significant ways like those with a negative Schwarzian derivative. It is surprising that this class of functions with an eventual negative Schwarzian derivative has not been previously considered as it naturally combines the seemingly nondynamic condition of having a negative Schwarzian derivative with the dynamics of the function. One result of this is that this enlarges the class of known functions with local and global properties similar to those of a function with a negative Schwarzian derivative. Moreover, as having an eventual negative Schwarzian derivative is not an asymptotic property, verification of this can often be done by direct computation. This last property makes functions with an eventual negative Schwarzian derivative a potentially useful tool in dealing with applications. Recall, that the introduction of this class was motivated by the maps arising in some models in neuroscience.

\section{Acknowledgments} The author would like to thank L. Bunimovich for suggesting the study of this problem as well as for useful discussions on the topic.

\end{document}